\documentclass[reqno]{amsart}
\usepackage{amssymb}
\usepackage{hyperref}
\usepackage{graphicx}
\usepackage{float}

\newtheorem{theorem}{Theorem}

\theoremstyle{remark}

\numberwithin{equation}{section}

\begin{document}

\title[Roots of continuous Hahn and Wilson polynomials]
{Stable equilibria for the roots of the symmetric continuous  Hahn and Wilson polynomials}

\author{J.F.  van Diejen}

\address{
Instituto de Matem\'atica y F\'{\i}sica, Universidad de Talca,
Casilla 747, Talca, Chile}

\email{diejen@inst-mat.utalca.cl}


\thanks{This work was supported in part by the {\em Fondo Nacional de Desarrollo
Cient\'{\i}fico y Tecnol\'ogico (FONDECYT)} Grant  \# 1170179.}

\date{May 2019}

\begin{abstract}
We show that the gradient flows associated with a  recently found family of Morse functions converge exponentially to the roots of the symmetric continuous Hahn polynomials. 
By symmetry reduction  the rate of the exponential convergence can be improved, which is clarified by comparing with corresponding gradient flows for the roots of
the Wilson polynomials.
\end{abstract}

\maketitle



\section{Introduction}\label{sec1}
For $\alpha,\beta >-1$, the Jacobi polynomials (cf. e.g. \cite[Chapter 9.8]{koe-les-swa:hypergeometric})
\begin{subequations}
\begin{equation}
P_n^{(\alpha ,\beta)} (x ):=\frac{ (\alpha+1)_n}{n!} \,
  {}_2F_1
\left[ \begin{array}{c} -n,n+ \alpha+\beta \\
\alpha +1 \end{array}  ; \frac{1-x}{2} \right] 
\end{equation}
satisfy the orthogonality relation
\begin{align}
\int_{-1}^1 &(1-x)^\alpha (1+x)^\beta P_n^{(\alpha ,\beta)} (x )P_m^{(\alpha ,\beta)} (x ) \\
&=\begin{cases}
\frac{2^{\alpha+\beta+1}\Gamma(n+\alpha +1)\Gamma (n+\beta+1) }
{(2n+\alpha+\beta+1)  \Gamma(n+\alpha+\beta +1) n!}  &\ \text{if}\ m=n, \\
0 &\ \text{if}\  m\neq n,
\end{cases} \nonumber
\end{align}
for $n,m=0,1,2,\ldots$.
\end{subequations}
Here $\Gamma (\cdot)$ refers to the gamma function, and we have adopted standard conventions for the Pochhammer symbol
\begin{equation*}
(a)_n := \prod_{0\leq k <n} (a+k), \quad (a_1,a_2,\ldots ,a_r)_n:=(a_1)_n (a_2)_n\cdots (a_r)_n ,
\end{equation*}
and for the hypergeometric series
\begin{equation*}
 {}_{r+1}F_r
\left[ \begin{array}{c} a_1,\ldots,a_{r+1} \\
b_1,\ldots ,b_r \end{array}  ; z  \right]  :=
\sum_{k\geq 0}  \frac{(a_1,\ldots,a_{r+1})_k}{(1,b_1,\ldots,b_r)_k}  z^k .
\end{equation*}
As an $n$th degree orthogonal polynomial, $P_n^{(\alpha ,\beta)} (x )$ has $n$ simple roots inside the interval of orthogonality:
\begin{equation*}
-1 < x_1^{(n)}<x_2^{(n)}<\cdots <x_n^{(n)}<1.
\end{equation*} 
In a landmark paper
\cite{sti:sur},  Stieltjes demonstrated that these roots $x_1^{(n)},\ldots ,x_n^{(n)}$
provide the coordinates of the
minimum of the potential
\begin{align}\label{electro}
V^{(\alpha,\beta)}(x_1,\ldots,x_n)= &- \sum_{1\leq j< k\leq n}   \log | x_j-x_k |  \\&    -\sum_{1\leq j\leq n}  \left( \frac{\alpha+1}{2} \log |x_j-1| + \frac{\beta+1}{2} \log |x_j+1| \right) \nonumber
\end{align}
in the domain $\{ (x_1,\ldots ,x_n) \in\mathbb{R}^n \mid -1<x_1<x_2<\cdots <x_n<1\}$.  
Since (by Gauss' law)
the strength of an electric field $E=-\nabla V$ caused by an infinite, straight, uniformly charged wire is proportional to the charge density and inversely proportional to the distance, the potential in Eq. \eqref{electro}
models the total energy of $n$ movable parallel wires with unit charge densities,
aligned perpendicular to the real axis between two additional  wires that are fixed at $x = 1$ and $x = -1$ with  positive charge densities $\frac{1}{2}(\alpha+1)$ and  $\frac{1}{2}(\beta+1)$, respectively.
Following Stieltjes' footsteps, electrostatic models of this type turn out to provide a very fruitful tool for the  study of the roots of classical orthogonal polynomials and their relatives, cf. e.g.
\cite{bel-mar-mar:some,dim-van:lame,for-rog:electrostatics,gru:variations,gru:electrostatic,hen-ros:electrostatic,hor:electrostatic,ism:electrostatics,ism:more,ism:classical,mar-mar-mar:electrostatic,sim:electrostatic,sze:orthogonal} and references therein.

As a member of the classical (hypergeometric) orthogonal families, the Jacobi polynomials satisfy a second-order differential equation \cite[Chapter 9.8]{koe-les-swa:hypergeometric}:
\begin{equation}
A(x) \frac{\text{d}^2}{\text{d}x^2}  P_n^{(\alpha ,\beta)} (x ) + B(x)  \frac{\text{d}}{\text{d}x}  P_n^{(\alpha ,\beta)} (x ) = \lambda_n    P_n^{(\alpha ,\beta)} (x ) ,
\end{equation}
where $A(x)=x^2-1$, $B(x)=(\alpha+1)(x+1)+(\beta+1)(x-1)$, and $\lambda_n=n(n+\alpha+\beta +1)$. In fact, such Sturm-Liouville type differential equations  play a fundamental role
in the proof of Stieltjes' electrostatic interpretation for the roots.
Moreover, in recent work of Steinerberger \cite{ste:electrostatic} it was pointed out that the  differential equations at issue also ensure
that the roots can be recovered via a dynamical system of the form:
\begin{align}\label{gs-jacobi}
\frac{\text{d}x_j}{\text{d}t}  =& - B(x_j)- A(x_j) \sum_{\substack{1\leq k\leq n\\ k\neq j}}     \frac{2}{x_j-x_k}  \\
=&  2 A(x_j)  \partial_{x_j} V^{(\alpha ,\beta)}(x_1,\ldots ,x_n)  
   \nonumber
\end{align}
($j=1,\ldots ,n$).  More specifically, as a special instance of  \cite[Theorem 2]{ste:electrostatic}, it follows that for any initial condition 
 $-1 < x_1(0)<\cdots <x_n(0)<1$ the solution $(x_1(t),\ldots ,x_n(t))$, $t\geq0$ of the system in Eq. \eqref{gs-jacobi} converges exponentially to the roots
 of the Jacobi polynomial $P_n^{(\alpha ,\beta)} (x ) $:
 \begin{equation}
 |x_j(t)-x_j^{(n)} |\leq K\exp (-\kappa_n t)  \qquad (j=1,\ldots ,n),
 \end{equation}
where $\kappa_n=\lambda_n-\lambda_{n-1}=2n+\alpha+\beta$ and $K$ denotes a positive constant whose value depends on  the parameters and on the initial condition.

Hypergeometric orthogonal polynomials beyond the Jacobi family generally do not fit within Stieltjes'  electrostatic framework, for instance because these often satisfy second-order difference equations rather than differential equations.
Still, for key examples of such families like the (symmetric) continuous Hahn and Wilson polynomials, it was recently observed that the roots can nevertheless be retrieved from the global minimum of a smooth Morse function,
which was found by drawing from ideas developed in the context of the Bethe Ansatz method for the diagonalization of quantum integrable particle models \cite{die-ems:solutions}.
In the present work it will be shown
that  the gradient flows associated with the pertinent Morse functions converge again exponentially to the roots in question. 
The underlying gradient system turns out to conserve parity symmetry (i.e., it preserves configurations that are invariant with respect to  coordinate reflection in the origin).
Upon performing the corresponding symmetry reduction in the case of the gradient system for the symmetric continuous Hahn polynomial, the rate of the convergence is further  improved, as becomes evident 
by comparing with the gradient system for  the roots of the Wilson polynomial.

To avoid confusion, it is imperative to recall at this point that the roots of the (basic) hypergeometric orthogonal polynomials belonging to the Askey scheme can also be characterized alternatively in terms of the equilibria 
of Calogero-Moser type integrable particle models  \cite{bih-cal:properties:a,bih-cal:properties:b,cal:equilibrium,die:equilibrium,oda-sas:equilibria,oda-sas:equilibrium,oda-sas:calogero,per:equilibrium}.
Specifically, the particle systems relevant for the roots of the continuous Hahn and Wilson families were introduced in Refs. \cite{die:deformations,die:difference}, and further analyzed at the level of quantum mechanics in \cite{die:multivariable} and at the level of classical mechanics in \cite{feh-gor:duality,pus:hyperbolic,pus:scattering}.
Here, however, we will not pursue this connection and rather try to stick closer to the original electrostatic interpretation of Stieltjes, which---for the polynomials of the symmetric continuous Hahn, Wilson and Askey-Wilson families---turns out to be governed by the Morse functions found in Ref. \cite{die-ems:solutions}.
A corresponding analysis of the gradient system for the roots of the  Askey-Wilson polynomial has  recently been carried out in Ref. \cite{die:gradient}.
A different approach to study the (extremal)  roots of the continuous Hahn polynomials and other ${}_3F_2$ reductions of the Wilson polynomials---based on contiguous relations derived via Zeilberger's algorithm---was presented in  Ref. \cite{joo-nji-koe:inner}.

In brief, the plan of the paper is as follows. The exponential convergence of the gradient flows for the roots of the symmetric continuous Hahn polynomials and the Wilson polynomials are established in Sections \ref{sec2} and \ref{sec3}, respectively.
In Section \ref{sec4} we conclude by pointing out how  the convergence rate in the case of the symmetric continuous Hahn polynomials can be improved by symmetry reduction.

\section{Symmetric Continuous Hahn Polynomials}\label{sec2}

\subsection{Preliminaries}
The  monic continuous Hahn polynomial $\text{p}_n(x ;a,b)$ \cite{ask-wil:set} of degree $n$ in $x$  
is,  for generic parameters, given by the terminating hypergeometric series
 \cite[\text{Chapter 9.4}]{koe-les-swa:hypergeometric}
\begin{align}\label{cHp}
\text{p}_n(x ;a,b)=&\frac{i^n (a+\overline{a},a+\overline{b})_n}{(n+a+\overline{a}+b+\overline{b}-1)_n} \\
 &\times  {}_3F_2
\left[ \begin{array}{c} -n,n+ a+\overline{a}+b+\overline{b}-1,a+ix  \\
a+\overline{a},a+\overline{b} \end{array}  ; 1 \right]  .\nonumber
\end{align}
For $\text{Re}(a),\text{Re}(b)>0$, these polynomials satisfy the orthogonality relations  \cite[\text{Chapter 9.4}]{koe-les-swa:hypergeometric}
\begin{subequations}
\begin{align}
&\frac{1}{2\pi} \int_{-\infty}^\infty  \Delta (x ;a,b)  \text{p}_n(x ;a,b) \text{p}_m(x ;a,b)  \text{d}x = \\
 &
\begin{cases}
\frac{n! \Gamma(a+\overline{a}+b+\overline{b}+ n-1)\Gamma ( a+	\overline{a}+n )\Gamma (a+\overline{b}+n) \Gamma (\overline{a}+b+n) \Gamma (b+\overline{b}+n) }
{\Gamma(a+\overline{a}+b+\overline{b}+2n-1) \Gamma (a+\overline{a}+b+\overline{b} +2n) }  &\ \text{if}\ m=n, \\
0 &\ \text{if}\  m\neq n,
\end{cases}
\nonumber
\end{align}
where
\begin{equation}\label{weight-cH}
\Delta (x;a,b)=\left| \Gamma (a+ix)\Gamma (b+ix) \right|^2.
\end{equation}
\end{subequations}
Upon assuming that possible non-real parameters  $a$ and $b$ arise as a complex conjugate pair, the weight function $\Delta (x;a,b)$ is even. The corresponding orthogonal polynomials are referred to as \emph{symmetric} continuous Hahn polynomials.

\subsection{Morse Function}\label{sub-morse-cH}
For parameter values belonging to the above orthogonality regime with possible non-real parameters  $a$ and $b$ arising as a complex conjugate pair, the roots
\begin{equation*}
x_1^{(n)}<x_2^{(n)}<\cdots <x_n^{(n)}
\end{equation*}
 of the symmetric continuous Hahn
polynomial $\text{p}_n(x ;a,b)$ turn out to minimize the Morse function \cite[Remarks 5.4, 5.5]{die-ems:solutions}
\begin{align}\label{VCH}
V(x_1,\ldots ,x_n;a,b) &= \sum_{1\leq j <k \leq n}   \int_0^{x_j-x_k} \arctan (\theta ) \text{d}\theta  \\
+ \sum_{1\leq j\leq n} &\left( \int_0^{x_j} \arctan \left( \frac{\theta}{a}\right)+  \arctan  \left( \frac{\theta}{b}\right)   \text{d}\theta +\frac{1}{2}\pi (n+1-2j )x_j \right) .
  \nonumber
\end{align}

To keep our presentation self-contained, let us briefly extract the main ingredients from the proof  in \cite{die-ems:solutions} sustaining the claim that the coordinates of  a critical point of
the above Morse function are necessarily  given by the roots of the symmetric continuous Hahn polynomial.
To this end we first observe that at such a critical point:
\begin{equation*}
{ \arctan\left(\frac{x_j}{a}\right)  +\arctan\left(\frac{x_j}{b}\right)  }  
+
 \sum_{\substack{1\leq k\leq n \\ k \neq j}} \arctan (x_j-x_k) ={ \frac{1}{2}\pi (2j-n-1) }
 \end{equation*}
 $(j=1,\ldots, n)$, which implies in particular that $x_1<x_2<\cdots <x_n$ in this situation. Moreover, 
since
\begin{equation*}
\exp \left(2i\arctan \bigl(x \bigr) \right)=\frac{i-x}{i+x}, 
\end{equation*}
it is also clear that then
\begin{equation}\label{bethe-cH}
\frac{\bigl(ia+x_j\bigr)\bigl(ib+x_j\bigr)}{\bigl(ia-x_j\bigr)\bigl(ib-x_j\bigr)}  \prod_{\substack{1\leq k \leq n\\ k\neq j}}  \frac{i+x_j-x_k}{i-x_j+x_k}=(-1)^{n+1}
\end{equation}
for $j=1,\ldots ,n$.
The upshot is that the asociated polynomial
\begin{equation*}
\text{p}_n(x)=(x-x_1)(x-x_2)\cdots (x-x_n)
\end{equation*}
satisfies the following second-order difference equation for the symmetric continuous Hahn polynomials \cite[\text{Chapter 9.4}]{koe-les-swa:hypergeometric}:
\begin{align}\label{de-cH}
& A(x;a,b) \Bigl( \text{p}_n(x+i) -\text{p}_n(x)\Bigr)+\\
&A(-x;a,b) \Bigl( \text{p}_n(x-i) -\text{p}_n(x)\Bigr) =\lambda_n  \text{p}_n(x), \nonumber
\end{align}
where $A(x;a,b)= (x+ia)(x+ib)$ and $\lambda_n=-n(n+2a+2b-1)$. Indeed, the fact that the basis of symmetric continuous Hahn polynomials $\text{p}_n(x;a,b)$, $n=0,1,2,\ldots$ satisfies
the difference equations in question confirms that the expression at the LHS of Eq. \eqref{de-cH} is actually a polynomial of degree $n$ in $x$ with a leading coefficient given by $\lambda_n$.
It thus suffices to establish the equality of both sides of the difference equation at the nodes $x_1,\ldots ,x_n$: 
\begin{equation*}\label{BAE-cH}
A  (x_j;a,b) \text{p}_n(x_j+i) +A(-x_j;a,b)\text{p}_n(x_j-i)=0\quad (j=1,\ldots ,n),
\end{equation*}
which amounts in turn to the relations in Eq. \eqref{bethe-cH}.
Hence $\text{p}_n(x)=\text{p}_n(x;a,b)$ (by the nondegeneracy of the eigenvalues $\lambda_n$ at the RHS of Eq. \eqref{de-cH}) and thus  $x_j=x_j^{(n)}$ ($j=1,\ldots ,n$).

It is manifest from the above argument  that the roots of the symmetric continuous Hahn polynomial $p_n(x;a,b)$ satisfy the identities in
\eqref{bethe-cH}, which is in fact an observation going back to 
\cite{oda-sas:equilibria}.

\subsection{Gradient Flow}
The corresponding gradient flow 
\begin{subequations}
\begin{equation}\label{g-flow}
\frac{\text{d}x_j}{\text{d}t} + \partial_{x_j} V(x_1,\ldots ,x_n;a,b)=0
\end{equation}
becomes of the form
\begin{align}\label{g-flow-chahn}
{  \frac{\text{d}x_j}{\text{d}t}+\frac{1}{2}\pi (n+1-2j) } &+
{ \arctan\left(\frac{x_j}{a}\right)  +\arctan\left(\frac{x_j}{b}\right)  }  \\
&+
 \sum_{\substack{1\leq k\leq n \\ k \neq j}} \arctan (x_j-x_k) =0 ,
 \nonumber
\end{align}
\end{subequations}
where $j=1,\ldots ,n$.
The following theorem affirms that the solutions of the gradient system \eqref{g-flow}, \eqref{g-flow-chahn}  converge exponentially to the equilibrium 
$\bigl( x_1^{(n)},x_2^{(n)}, \ldots , x_n^{(n)}\bigr)$ given by the symmetric continuous Hahn roots.

\begin{theorem}\label{cHzeros:thm}
\begin{subequations}
Let $\text{Re}(a),\text{Re}(b)>0$ with possible non-real parameters  $a$ and $b$ arising as a complex conjugate pair.
\begin{itemize}
\item[a)]  The unique global minimum in $\mathbb{R}^n$ of the strictly convex, radially unbounded,
Morse function $V(x_1,\ldots ,x_n;a,b)$ \eqref{VCH}  is attained at the symmetric continuous Hahn roots $x_j= x_j^{(n)}$ ($j=1,\ldots ,n$).
\item[b)]  Let
$\bigl(x_1(t),\ldots ,x_n(t)\bigr)$, $t\geq 0$ denote the unique solution of the gradient system \eqref{g-flow}, \eqref{g-flow-chahn} determined by a choice for
the initial condition 
\begin{equation*}
(x_1(0),\ldots ,x_n(0))\in\mathbb{R}^n,
\end{equation*}
 and let us fix any $\kappa$ in the interval
\begin{equation}
{ 0< \kappa <     \frac{\text{Re}(a)}{\text{Re}^2(a)+(R_n+|\text{Im}(a)|)^2}+ \frac{\text{Re}(b)}{\text{Re}^2(b)+(R_n + |\text{Im}(b)| )^2} },
\end{equation}
where $R_n:=x_n^{(n)}$. Then there exists a constant  $K>0$ (depending on the parameter values, on the initial condition, and on $\kappa$) such that 
\begin{equation}
|x_j(t)-x_j^{(n)}|\leq  K e^{-\kappa t} \qquad (j=1,\ldots,n)
\end{equation}
whenever $t\geq 0$.
\end{itemize}
\end{subequations}
\end{theorem}

\subsection{Proof of Theorem \ref{cHzeros:thm}}
Because $V(x_1,\ldots ,x_n;a,b)$ \eqref{VCH} is manifestly smooth on $\mathbb{R}^n$, the existence of a global minimum is immediate from the observation that the function under consideration is radially unbounded:
 $V(x_1,\ldots ,x_n;a,b)\to +\infty$ when $x_1^2+\cdots +x_n^2\to + \infty$ (cf. \cite[Section 5.1]{die-ems:solutions} for a detailed verification of this crucial property). The convexity reconfirms the uniqueness of the global minimum. Indeed,
 the relevant Hessian is of the form
\begin{align*}
&H_{j,k}(x_1,\ldots ,x_n;a,b) =\partial_{x_j}\partial_{x_{k}}V(x_1,\ldots ,x_n;a,b) \\
&=
\begin{cases}
 \frac{a}{a^2+x_j^2}+ \frac{b}{b^2+x_j^2}+
 \sum_{l \neq j}  \frac{1}{1+(x_j-x_l)^2} & \text{if $j=k$,}\\
   -\frac{1}{1+(x_j-x_k)^2}  & \text{if $j\neq k$,}\\
\end{cases} 
\nonumber
\end{align*}
so
\begin{align*}
\sum_{1\leq j,k \leq n}  y_j y_{k}H_{j,k} (x_1,\ldots,x_n;a,b)  
=& \sum_{1\leq j\leq n} \left(\frac{a}{a^2+x_j^2}+ \frac{b}{b^2+x_j^2}
 \right) y_j^2  \nonumber \\
&+  \sum_{1\leq j< k \leq n } 
\frac{(y_j-y_k)^2 }{1+  (x_j-x_k)^2} >0  
\end{align*}
when $y_1^2+\cdots +y_n^2\neq 0$. The fact that this unique global minimum is attained precisely at the roots of the symmetric continuous Hahn polynomial $\text{p}_n(x;a,b)$ follows via
\cite[Remarks 5.4, 5.5]{die-ems:solutions} as outlined above in Subsection \ref{sub-morse-cH}, which completes the proof of  part a) of the theorem.

We conclude from part a) that $(x_1^{(n)},\ldots,x_n^{(n)})$ constitutes the unique equilibrium of our gradient system  \eqref{g-flow}, \eqref{g-flow-chahn}. 
Given any initial condition
in $\mathbb{R}^n$, the smoothness of the Morse function guarantees 
 the existence and uniqueness of the solution $(x_1(t),\ldots,x_n(t))$ locally for small $t\geq 0$.
Moreover, since $V(x_1,\ldots, x_n;a,b)$ is radially unbounded and strictly decreasing along the gradient flow:
\begin{equation*}
\frac{\text{d}}{\text{d}t} V \bigl (x_1(t), \ldots ,x_n(t);a,b\bigr)= -\sum_{j=1}^n \left| (\partial_{x_j} V) \bigl (x_1(t), \ldots ,x_n(t);a,b\bigr) \right|^2 < 0 
\end{equation*}
(outside the equilibrium),  the solution cannot escape to infinity in finite time and thus extends (uniquely) for all $t\geq 0$. In fact, the radial unboundedness of the
potential $V(x_1,\ldots, x_n;a,b)$ \eqref{VCH} guarantees that the equilibrium $(x_1^{(n)},\ldots,x_n^{(n)})$ is globally asymptotically stable (cf. e.g. \cite[Theorem 4.2]{kha:nonlinear}):
\begin{equation}\label{asymptotic}
\lim_{t\to +\infty}   x_j(t) = x_j^{(n)}\qquad (j=1,\ldots ,n) .
\end{equation}
The rate of the convergence to the equilibrium is determined by the lowest eigenvalue $\lambda$ of the Jacobian of the linearized system at the equilibrium point. This Jacobian is given by
the above Hessian.  Since the weight function $\Delta (x;a,b)$  \eqref{weight-cH} is even in $x$, the roots of $\text{p}_n(x;a,b)$ are symmetrically distributed around the origin. Hence
 $-R_n\leq x_j\leq R_n$ ($j=1,\ldots ,n$) at the equilibrium, and thus
 \begin{align*}
  \text{Re} \Bigl( \frac{a}{a^2+x_j^2} \Bigr) =&
\frac{1}{2}\left( \frac{\text{Re}(a)}{\text{Re}^2(a)+(x_j+\text{Im}(a))^2}+ \frac{\text{Re}(a)}{\text{Re}^2(a)+(x_j-\text{Im}(a))^2}\right) \\
\geq &  \frac{\text{Re}(a)}{\text{Re}^2(a)+(R_n+|\text{Im}(a)|)^2}
 \end{align*}
in this situation. The upshot is that
\begin{align*}
&\sum_{1\leq j,k \leq n}  y_j y_{k}H_{j,k} (x_1^{(n)},\ldots,x_n^{(n)};a,b)    \\ 
& \geq {\textstyle \left( \frac{\text{Re}(a)}{\text{Re}^2(a)+(R_n+|\text{Im}(a)|)^2}+ \frac{\text{Re}(b)}{\text{Re}^2(b)+(R_n+|\text{Im}(b)|)^2}\right) ( y_1^2+\cdots +y_n^2) }+  \sum_{1\leq j< k \leq n } 
\frac{(y_j-y_k)^2 }{1+  4R_n^2 }   \\
&\geq  {\textstyle \left( \frac{\text{Re}(a)}{\text{Re}^2(a)+(R_n+|\text{Im}(a)|)^2}+ \frac{\text{Re}(b)}{\text{Re}^2(b)+(R_n+|\text{Im}(b)|)^2}\right) ( y_1^2+\cdots +y_n^2) },
\end{align*}
which entails that
\begin{align*}
\lambda &= \min_{\substack{(y_1,\ldots,y_n)\in\mathbb{R}^n\\ y_1^2+\cdots +y_n^2\neq 0}}    \frac{\sum_{1\leq j,k\leq n} y_jy_k H_{j,k}(x_1^{(n)},\ldots ,x_n^{(n)};a,b)}{y_1^2+\cdots +y_n^2}\\
&\geq 
 \frac{\text{Re}(a)}{\text{Re}^2(a)+(R+|\text{Im}(a)|)^2}+ \frac{\text{Re}(b)}{\text{Re}^2(b)+(R + |\text{Im}(b)| )^2} > \kappa .
\end{align*}
This lower bound on the eigenvalue guarantees  (cf. e.g. \cite[Corollary 2.78]{chi:ordinary})  that there exists a neighborhood $U\subset\mathbb{R}^n$ around the equilibrium point $\bigl( x_1^{(n)}, \ldots , x_n^{(n)}\bigr)$ and a constant $C>0$ such that
\begin{equation*}
|x_j(t)-x_j^{(n)}|\leq C |x_j(t_0)-x_j^{(n)}|  e^{-\kappa (t-t_0)}   \qquad (j=1,\ldots,n),
\end{equation*}
whenever
$(x_1(t_0),\ldots ,x_n(t_0))\in U$ and
 $t\geq t_0$. To complete the proof of part b), it now suffices to recall that the gradient flow pulls us from any initial condition $(x_1(0),\ldots ,x_n(0))\in \mathbb{R}^n$ into $U$ within finite time $t$ (in view of  Eq. \eqref{asymptotic}).

\subsection{Numerical Samples}\label{cH-numerics}
\begin{figure}[h]
 \centering
 \caption{Continuous Hahn trajectories $x^{(30)}_1(t),\ldots, x_{30}^{(30)}(t)$   corresponding to the initial condition $x_j^{(30)}(0)=0$ ($j=1,\ldots ,30$) and the parameter values $a=10$ and $b=\frac{3}{10}$.}
\includegraphics[scale=0.7]{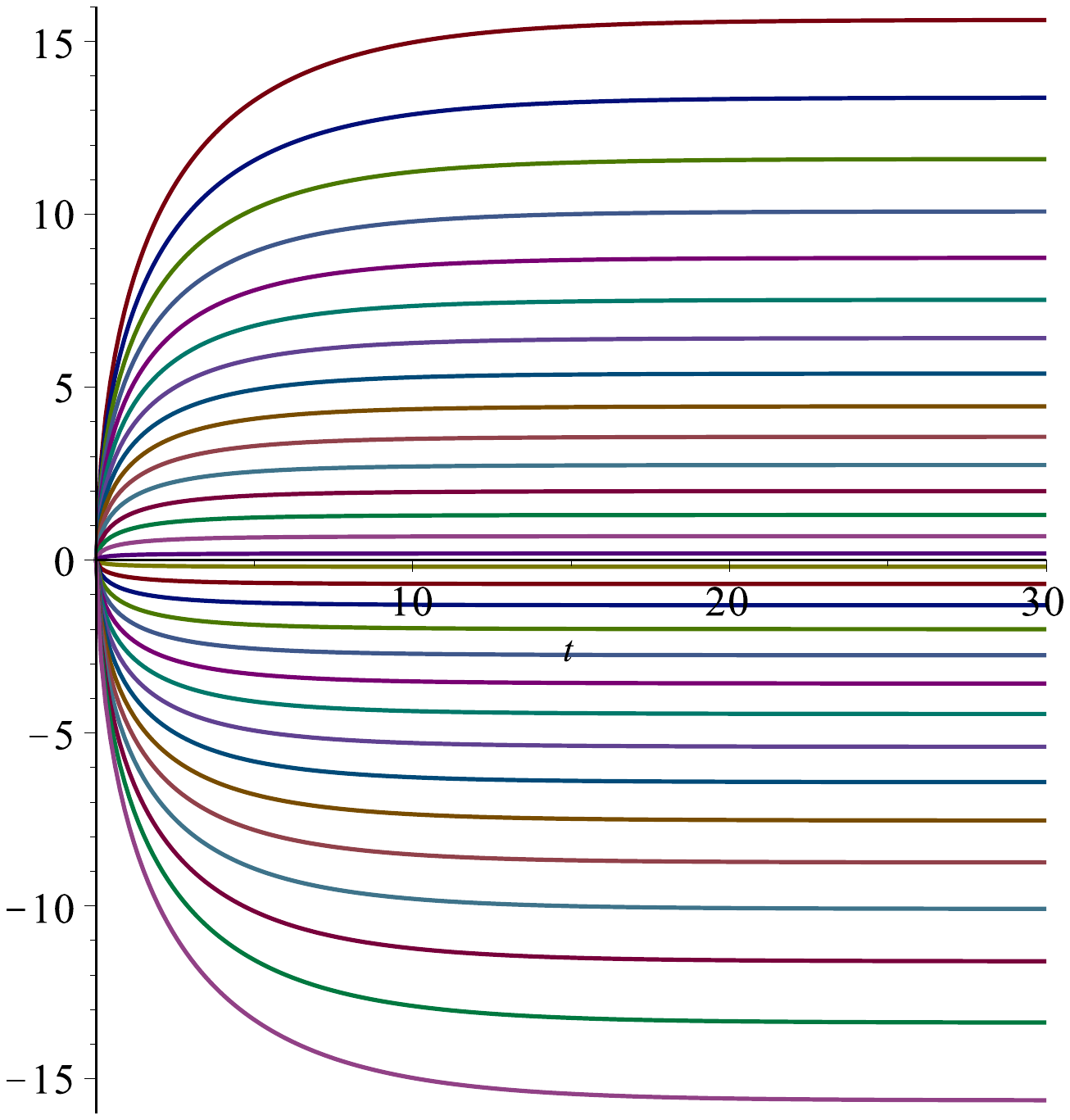} \label{cH-fig1}
\end{figure}

\begin{figure}[h]
  \centering
  \caption{Evolution of the continuous Hahn logarithmic error  $\log\bigl |x_j(t)-x_j^{(30)}\bigr| $ for $j=1$ (top), $j=8$ (middle) and $j=15$ (bottom),
  corresponding to the initial condition $x_j^{(30)}(0)=0$ ($j=1,\ldots ,30$) and the parameter values $a=10$ and $b=\frac{3}{10}$.}
\includegraphics[scale=0.7]{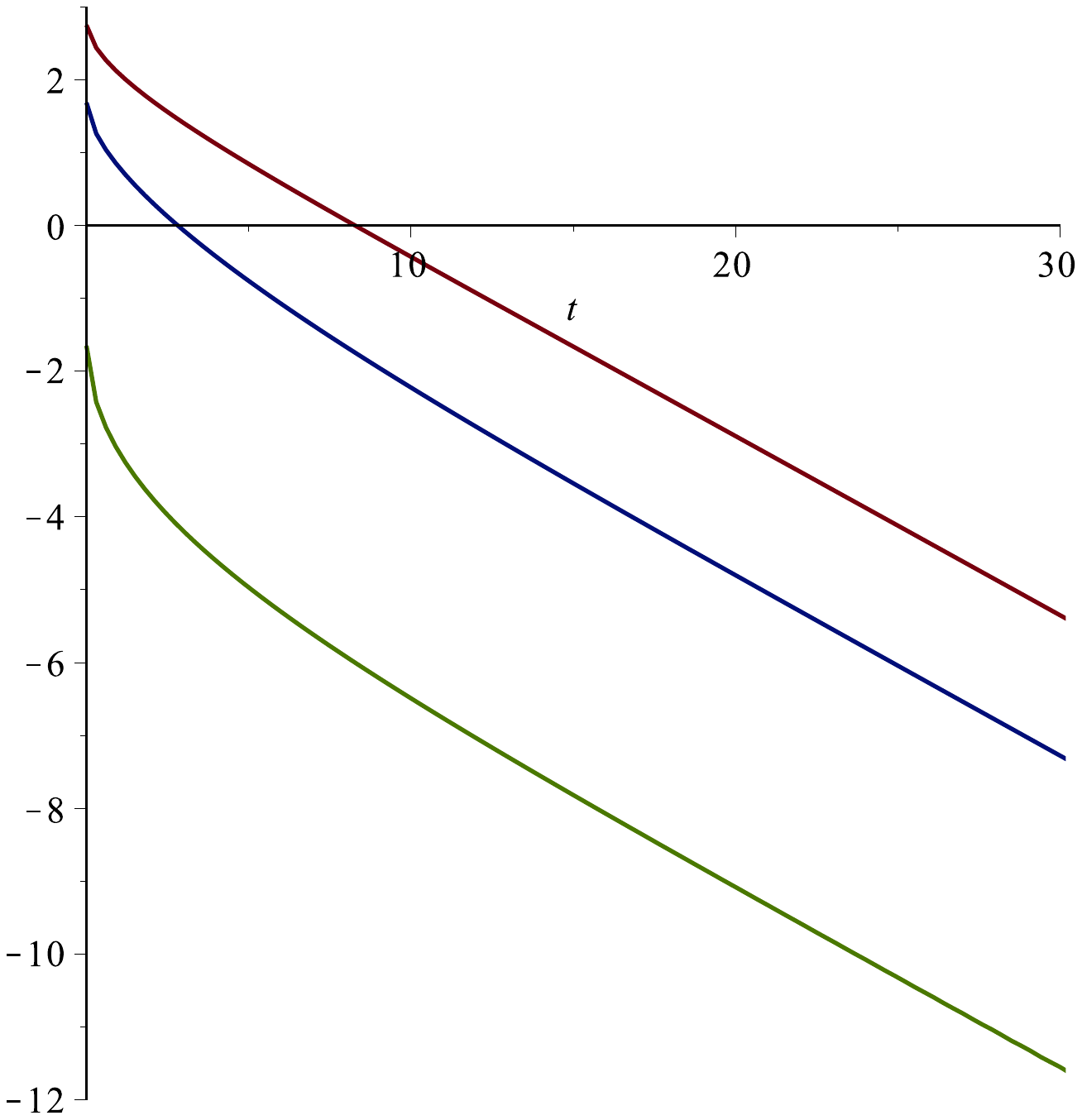}  \label{cH-fig2}
\end{figure}

For $n=30$ with
$a=10$ and $b=\frac{3}{10}$, the continuous Hahn roots  become with a precision of 4 decimals:
\begin{equation*}
\begin{matrix}
x_1^{(30)}=-15.6230 ,& & x_6^{(30)}=-7.5285 ,     && x_{11}^{(30)}=-2.7503 ,  \\
 x_2^{(30)}=-13.3738 ,& &  x_7^{(30)}=-6.4188 , &&  x_{12}^{(30)}=-1.9957 , \\
 x_3^{(30)}=-11.6001, & &  x_8^{(30)}=-5.3956 , && x_{13}^{(30)}=-1.3059 ,  \\
  x_4^{(30)}=-10.0841 ,& &   x_9^{(30)}=-4.4474 ,&& x_{14}^{(30)}=-0.6907 ,  \\
 x_5^{(30)}=-8.7415 ,& &    x_{10}^{(30)}=-3.5671 , &&  x_{15}^{(30)}=-0.1919     
 \end{matrix}
\end{equation*}
(and $x^{(30)}_{j}=-x^{(30)}_{31-j}$ for $j=16,\ldots ,30$).
The corresponding trajectories of $x_j^{(n)}(t)$ for $0\leq t\leq 30$ starting from the origin, i.e. with $x_j^{(n)}(0)=0$ ($j=1,\ldots ,n$), are exhibited in Figure \ref{cH-fig1}.
The slope of the evolution of the logarithmic error $\log \bigl| x_j(t)-x_j^{(n)}\bigr|$ in Figure \ref{cH-fig2} confirms that the convergence is exponential with a decay rate that
exceeds the not very sharp estimate of $\frac{a}{a^2+R_n^2}+\frac{b}{b^2+R_n^2}\approx 0.030$ 
guaranteed by Theorem \ref{cHzeros:thm}.

\section{Wilson Polynomials}\label{sec3}

\subsection{Preliminaries}
The  monic Wilson polynomial $\text{p}_n(x^2 ;a,b,c,d)$  \cite{wil:some} is an orthogonal  polynomial of degree $n$ in $x^2$ that sits at the top of Askey's hypergeometric scheme \cite{koe-les-swa:hypergeometric}.
For generic parameter  values it is given by the terminating hypergeometric series
\cite[\text{Chapter 9.1}]{koe-les-swa:hypergeometric}
\begin{align}\label{Wp}
\text{p}_n(x^2 ;a,b,c,d)=&\frac{(-1)^n (a+b,a+c,a+d)_n}{(n+a+b+c+d-1)_n} \\
 &\times  {}_4F_3
\left[ \begin{array}{c} -n,n+a+b+c+d-1,a+ix ,a-ix \\
a+b,a+c,a+d \end{array}  ; 1 \right] ,\nonumber
\end{align}
When $\text{Re}(a), \text{Re}(b),\text{Re}(c),\text{Re}(d)>0$ with possible non-real parameters arising in complex conjugate pairs, the polynomials in question satisfy the following orthogonality relations
\cite[\text{Chapter 9.1}]{koe-les-swa:hypergeometric}
\begin{subequations}
\begin{align}
&\frac{1}{2\pi} \int_0^\infty  \Delta (x ;a,b,c,d)  \text{p}_n(x^2 ;a,b,c,d) \text{p}_m(x ^2;a,b,c,d)  \text{d}x = \\
 &
\begin{cases}
\frac{n! \Gamma(a+b+c+d+ n-1)\Gamma ( a+b+n )\Gamma (a+c+n) \Gamma (a+d+n) \Gamma (b+c+n)  \Gamma (b+d+n) \Gamma (c+d+n) }
{\Gamma(a+b+c+d+2n-1) \Gamma (a+b+c+d +2n) }  &\ \text{if}\ m=n, \\
0 &\ \text{if}\  m\neq n,
\end{cases}
\nonumber
\end{align}
where
\begin{equation}\label{weight-W}
\Delta (x ;a,b,c,d)=\left| \frac{\Gamma (a+ix)\Gamma (b+ix)\Gamma (c+ix)\Gamma (d+ix)} {\Gamma (2ix )}\right|^2 .
\end{equation}
\end{subequations}

\subsection{Morse Function}\label{sub-morse-W}
For parameters values within the above orthogonality regime, the roots
\begin{equation*}
0<x_1^{(n)}<x_2^{(n)}<\cdots <x_n^{(n)}
\end{equation*}
of the Wilson polynomial $\text{p}_n(x^2 ;a,b,c,d)$ 
minimize the Morse function \cite[Remarks 5.4, 5.5]{die-ems:solutions}
\begin{align}
V(x_1,\ldots ,x_n;a,b,c,d) =  \sum_{1\leq j <k \leq n} & \left( \int_0^{x_j+x_k} \arctan (\theta ) \text{d}\theta + \int_0^{x_j-x_k} \arctan (\theta ) \text{d}\theta \right)  \nonumber \\
+ \sum_{1\leq j\leq n} &\left( \int_0^{x_j}  \arctan \left( \frac{\theta}{a}\right) + 
 \arctan  \left( \frac{\theta}{b}\right)  \right.  \label{VW}  \\
&\qquad  \left.  + \arctan \left( \frac{\theta}{c}\right) 
+ \arctan \left( \frac{\theta}{d}\right) \text{d}\theta - \pi j x_j \right) .
  \nonumber
\end{align}

Indeed, we now have that at the critical points of the Morse function
\begin{align*}
\arctan\left(\frac{x_j}{a}\right) & +\arctan\left(\frac{x_j}{b}\right) +\arctan\left(\frac{x_j}{c}\right)+\arctan\left(\frac{x_j}{d}\right)  \\
&+
 \sum_{\substack{1\leq k\leq n \\ k \neq j}} \bigl( \arctan (x_j+x_k)+\arctan (x_j-x_k) \bigr) =\pi j
 \nonumber
\end{align*}
($j=1,\ldots ,n$), so $0<x_1<x_2<\cdots <x_n$ and (after exponentiation)
\begin{equation*}\label{bethe-W}
\frac{\bigl(ia+x_j\bigr)\bigl(ib+x_j\bigr)\bigl(ic+x_j\bigr)\bigl(id+x_j\bigr)}{\bigl(ia-x_j\bigr)\bigl(ib-x_k\bigr)\bigl(ic-x_j\bigr)\bigl(id-x_j\bigr)}  \prod_{\substack{1\leq k \leq n\\ k\neq j}}  \frac{(i+x_j+x_k)(i+x_j-x_k)}{(i-x_j-x_k)(i-x_j+x_k)}=1
\end{equation*}
for $j=1,\ldots ,n$. Upon rewriting the latter identities in the form
\begin{equation*}
A  (x_j;a,b,c,d) \text{p}_n(x_j+i) +A(-x_j;a,b,c,d)\text{p}_n(x_j-i)=0\quad (j=1,\ldots ,n),
\end{equation*}
with $A(x;a,b,c,d)= \frac{(x+ia)(x+ib)(x+ic)(x+id)}{2x(2x+i)}$ and 
\begin{equation*}
\text{p}_n(x)=(x^2-x^2_1)(x^2-x^2_2)\cdots (x^2-x^2_n),
\end{equation*}
 it is deduced
in the same way as before that the factorized polynomial in question
satisfies the second-order difference equation for the Wilson polynomials \cite[\text{Chapter 9.1}]{koe-les-swa:hypergeometric}:
\begin{align*}
& A(x;a,b,c,d) \Bigl( \text{p}_n(x+i) -\text{p}_n(x)\Bigr)+\\
&A(-x;a,b,c,d) \Bigl( \text{p}_n(x-i) -\text{p}_n(x)\Bigr) =\lambda_n \text{p}_n(x) \nonumber
\end{align*}
with $\lambda_n=-n(n+a+b+c+d-1) $.
 The upshot is that $\text{p}_n(x)=\text{p}_n(x^2;a,b,c,d)$ (using again the nondegeneracy of the eigenvalues $\lambda_n$ at the RHS) and thus  $x_j=x_j^{(n)}$ ($j=1,\ldots ,n$).

Let us stress in this connection that the observation that the roots of Wilson polynomial $p_n(x^2;a,b,c,d)$ satisfy the algebraic identities in
the second equation below the Morse function \eqref{VW} goes again back to 
\cite{oda-sas:equilibria} (cf. also \cite{die:equilibrium} and \cite{bih-cal:properties:a}).

\subsection{Gradient Flow}
The corresponding gradient flow now becomes:
\begin{align}\label{g-flow-wilson}
{ \frac{\text{d} x_j}{\text{d} t}-\pi j } &+
{ \arctan\left(\frac{x_j}{a}\right)  +\arctan\left(\frac{x_j}{b}\right) +\arctan\left(\frac{x_j}{c}\right)+\arctan\left(\frac{x_j}{d}\right)} \\
&+
 \sum_{\substack{1\leq k\leq n \\ k \neq j}} \bigl( \arctan (x_j+x_k)+\arctan (x_j-x_k) \bigr) =0 ,
 \nonumber
\end{align}
$j=1,\ldots ,n$.
The following theorem affirms that the solutions of the gradient system \eqref{g-flow-wilson} converge exponentially to the equilibrium 
$\bigl( x_1^{(n)}, \ldots , x_n^{(n)}\bigr)$ given by the Wilson roots.

\begin{theorem}\label{Wzeros:thm}
\begin{subequations}
Let $\text{Re}(a),\text{Re}(b),\text{Re}(c),\text{Re}(d)>0$  with possible non-real parameters arising in complex conjugate pairs.
\begin{itemize}
\item[a)]  The unique global minimum of the strictly convex, radially unbounded,
Morse function $V(x_1,\ldots ,x_n;a,b,c,d)$ \eqref{VW}  is attained at the Wilson roots $x_j= x_j^{(n)}$ ($j=1,\ldots ,n$).
\item[b)]  Let
$\bigl(x_1(t),\ldots ,x_n(t)\bigr)$, $t\geq 0$ denote the unique solution of the gradient system \eqref{g-flow-wilson} determined by a choice for
the initial condition
\begin{equation*}
(x_1(0),\ldots ,x_n(0))\in\mathbb{R}^n, 
\end{equation*}
and let us fix any $\kappa$ in the interval
\begin{align}\label{W-rate}
{\textstyle 0< \kappa <   \frac{2(n-1)}{1+4R_n^2}+  \frac{\text{Re}(a)}{\text{Re}^2(a)+(R_n+|\text{Im}(a)|)^2}+ \frac{\text{Re}(b)}{\text{Re}^2(b)+(R_n+|\text{Im}(b)|)^2} } & \\
{\textstyle + \frac{\text{Re}(c)}{\text{Re}^2(c)+(R_n+|\text{Im}(c)|)^2}+ \frac{\text{Re}(d)}{\text{Re}^2(d)+(R_n+|\text{Im}(d)|)^2} }& \nonumber
\end{align}
where $R_n:=x_n^{(n)}$. Then there exists a constant  $K>0$ (depending on the parameter values, on the initial condition, and on $\kappa$) such that 
\begin{equation}
|x_j(t)-x_j^{(n)}|\leq  K e^{-\kappa t} \qquad (j=1,\ldots,n)
\end{equation}
whenever $t\geq 0$.
\end{itemize}
\end{subequations}
\end{theorem}

\subsection{Proof of Theorem \ref{Wzeros:thm}}
The proof runs along the same lines as that of Theorem \ref{cHzeros:thm}, so we only highlight some of the principal modifications in the corresponding computations.

As before, the first part of the theorem hinges on  the considerations in
\cite[Remarks 5.4, 5.5]{die-ems:solutions}  with the main points summarized above in Subsection \ref{sub-morse-W}
(cf. also the proof of \cite[Proposition 4.1]{die-ems:solutions} for a detailed check that the present Morse function is indeed radially unbounded).

To infer the asserted estimate for the convergence rate, we must again provide the lower bound for the eigenvalues of the pertinent
Jacobian evaluated at the equilibrium point.  This Jacobian is given by the Hessian
\begin{align}\label{Hesse-AW}
&H_{j,k}(x_1,\ldots ,x_n;a,b,c,d) =\partial_{x_j}\partial_{x_{k}}V(x_1,\ldots ,x_n;a,b,c,d) \\
&=
\begin{cases}
\frac{a}{a^2+x_j^2}  +\frac{b}{b^2+x_j^2} +\frac{c}{c^2+x_j^2}  +\frac{d}{d^2+x_j^2}  +
 \sum_{l \neq j} \bigl( \frac{1}{1+(x_j+x_l)^2}+ \frac{1}{1+(x_j-x_l)^2} \bigr) & \text{if $j=k$,}\\
 \frac{1}{1+(x_j+x_l)^2} - \frac{1}{1+(x_j-x_l)^2} & \text{if $j\neq k$,}\\
\end{cases} 
\nonumber
\end{align}
which confirms the convexity:
\begin{align*}
\sum_{1\leq j,k \leq n}  y_j y_{k} H_{j,k} (x_1,\ldots,x_n;a,b,c,d)  
= \sum_{1\leq j\leq n} \Bigl(  \sum_{\epsilon\in \{a,b,c,d\}}  \frac{\epsilon}{\epsilon^2+x_j^2}
 \Bigr) y_j^2 &  \nonumber \\
+  \sum_{1\leq j< k \leq n } \Bigl(
\frac{(y_j+y_k)^2 }{1+  (x_j+x_k)^2} +\frac{(y_j-y_k)^2 }{1+  (x_j-x_k)^2} \Bigr) . &
\end{align*}
Since $0<x_j\leq R_n$ ($j=1,\ldots ,n$) at the equilibrium, we now see that
\begin{align*}
&\sum_{1\leq j,k \leq n}  y_j y_{k}H_{j,k} (x_1^{(n)},\ldots,x_n^{(n)};a,b,c,d)    \\ 
& \geq {\textstyle \left(  \sum_{\epsilon\in \{a,b,c,d\}} \frac{\text{Re}(\epsilon)}{\text{Re}^2(\epsilon)+(R_n+|\text{Im}(\epsilon)|)^2}\right) ( y_1^2+\cdots +y_n^2) } \\
&+  \sum_{1\leq j< k \leq n } 
\frac{(y_j+y_k)^2+(y_j-y_k)^2 }{1+  4R_n^2 }   \\
&=  {\textstyle \left(    \sum_{\epsilon\in \{a,b,c,d\}} \frac{\text{Re}(\epsilon)}{\text{Re}^2(\epsilon)+(R_n+|\text{Im}(\epsilon)|)^2}  +\frac{2(n-1)}{1+4R_n^2} \right) ( y_1^2+\cdots +y_n^2) }.
\end{align*}
This entails the desired lower bound for the smallest eigenvalue:
\begin{align*}
\lambda &= \min_{\substack{(y_1,\ldots,y_n)\in\mathbb{R}^n\\ y_1^2+\cdots +y_n^2\neq 0}}    \frac{\sum_{1\leq j,k\leq n} y_jy_k H_{j,k}(x_1^{(n)},\ldots ,x_n^{(n)};a,b,c,d)}{y_1^2+\cdots +y_n^2}\\
&\geq 
   \sum_{\epsilon\in \{a,b,c,d\}} \frac{\text{Re}(\epsilon)}{\text{Re}^2(\epsilon)+(R_n+|\text{Im}(\epsilon)|)^2}  +\frac{2(n-1)}{1+4R_n^2}> \kappa .
\end{align*}

\subsection{Numerical Samples}
\begin{figure}[h]
 \centering
 \caption{Wilson trajectories $x^{(15)}_1(t),\ldots, x_{15}^{(15)}(t)$  corresponding to the  initial condition $x_j^{(15)}(0)=0$ ($j=1,\ldots ,15$) and the parameter values
 $a=\frac{17}{3}$, $b=\frac{1}{5}$, $c=1+i$, and $d=1-i$.}
\includegraphics[scale=0.7]{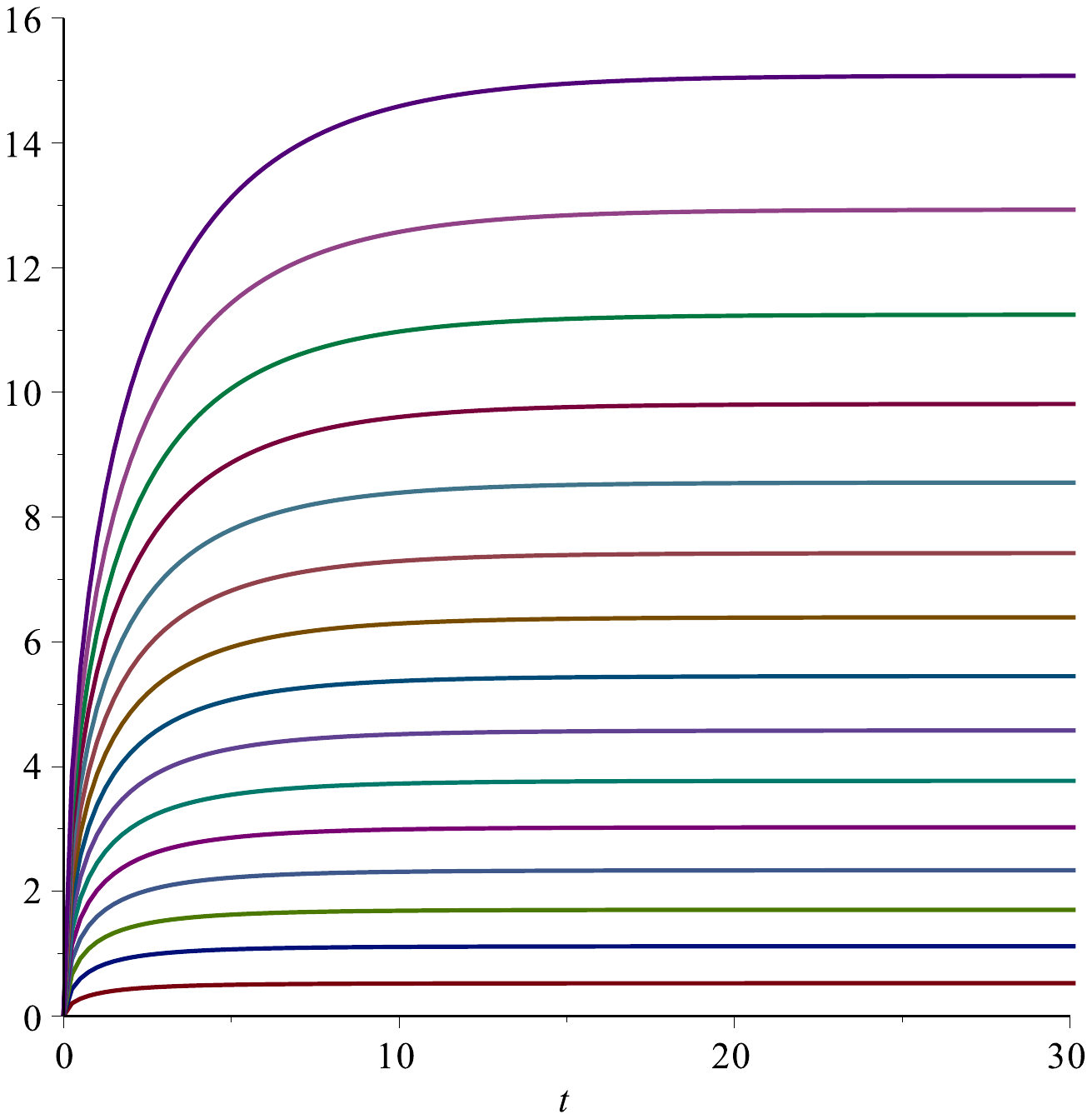} \label{W-fig1}
\end{figure}

\begin{figure}[h]
  \centering
 \caption{Evolution of the Wilson logarithmic error  $\log\bigl |x_j(t)-x_j^{(15)}\bigr| $ for $j=1$ (bottom), $j=8$ (middle) and $j=15$ (top), corresponding to the  initial condition $x_j^{(15)}(0)=0$ ($j=1,\ldots ,15$) and the parameter values
 $a=\frac{17}{3}$, $b=\frac{1}{5}$, $c=1+i$, and $d=1-i$.}
\includegraphics[scale=0.7]{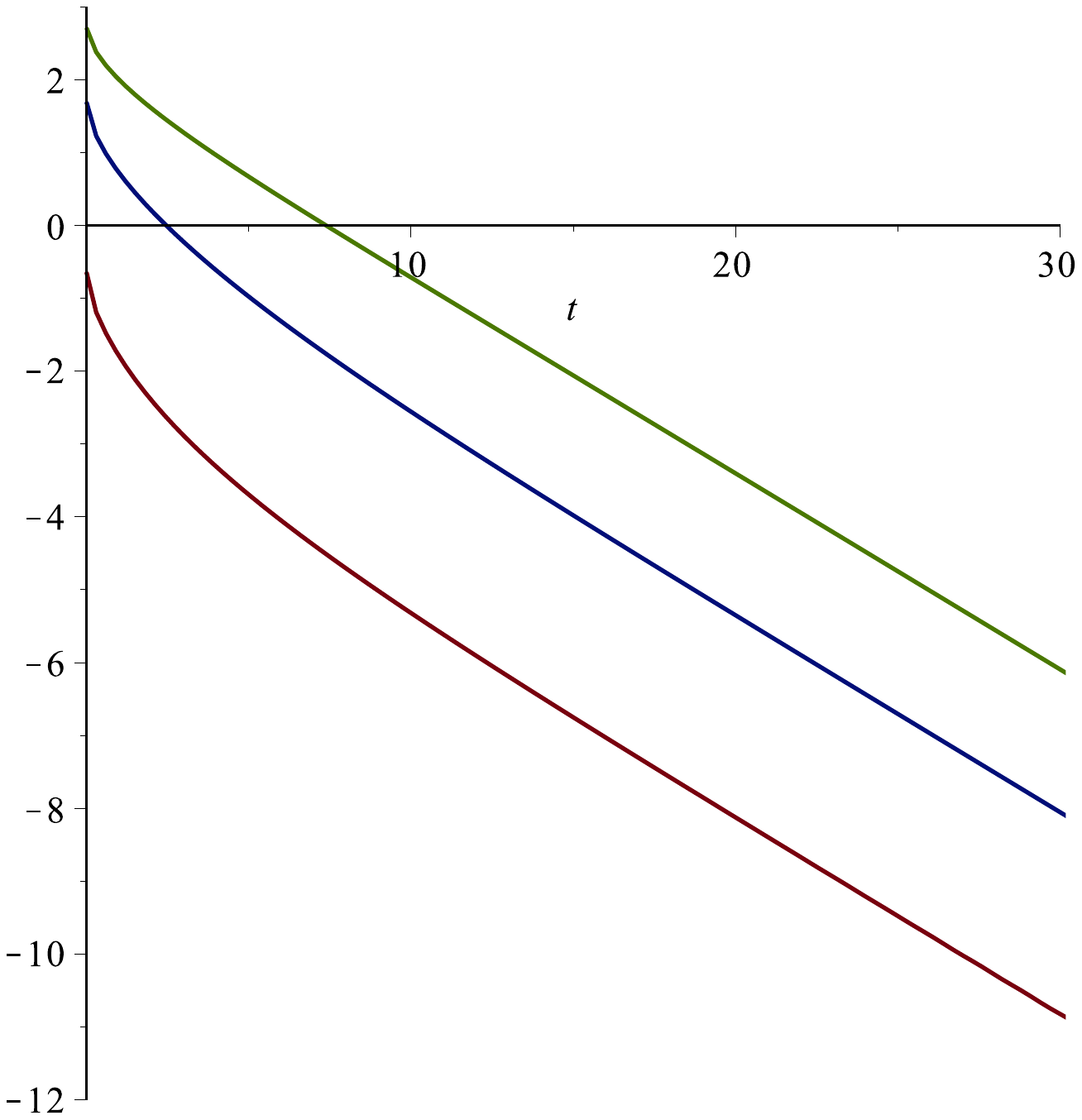}  \label{W-fig2}
\end{figure}

For $n=15$ with
$a=\frac{17}{3}$, $b=\frac{1}{5}$, $c=1+i$, and $d=1-i$, the Wilson roots  become with a precision of 4 decimals:
\begin{equation*}
\begin{matrix}
x_1^{(15)}=0.5274 ,& & x_6^{(15)}=3.7728 ,     && x_{11}^{(15)}=8.5546 ,  \\
 x_2^{(15)}=1.1194 ,& &  x_7^{(15)}=4.5787 , &&  x_{12}^{(15)}=9.8143 , \\
 x_3^{(15)}=1.7050, & &  x_8^{(15)}=5.4496 , && x_{13}^{(15)}=11.2449 ,  \\
  x_4^{(15)}=2.3375 ,& &   x_9^{(15)}=6.3938 ,&& x_{14}^{(15)}=12.9284 ,  \\
 x_5^{(15)}=3.0266 ,& &    x_{10}^{(15)}=7.4231 , &&  x_{15}^{(15)}=15.0759  .  
 \end{matrix}
\end{equation*}
The corresponding trajectories of $x_j^{(n)}(t)$ for $0\leq t\leq 30$, with an initial condition  of the form $x_j^{(n)}(0)=0$ ($j=1,\ldots ,n$), are exhibited in Figure \ref{W-fig1}.
The slope of the evolution of the logarithmic error $\log \bigl| x_j(t)-x_j^{(n)}\bigr|$ in Figure \ref{W-fig2} confirms that the convergence is exponential with a decay rate that
exceeds the not very sharp estimate of
\begin{equation*}
 \sum_{\epsilon\in \{a,b,c,d\}} \frac{\text{Re}(\epsilon)}{\text{Re}^2(\epsilon)+(R_n+|\text{Im}(\epsilon)|)^2}  +\frac{2(n-1)}{1+4R_n^2} \approx 0.061
\end{equation*}
guaranteed by Theorem \ref{Wzeros:thm}.

\section{Symmetry Reduction}\label{sec4}
As noticed above, the weight function $\Delta (x;a,b)$ \eqref{weight-cH} is even in $x$ for our parameter regime, so the roots of the corresponding symmetric continuous Hahn polynomial $\text{p}_n(x;a,b)$ are symmetrically distributed around the origin:
\begin{equation}
x_j^{(n)}+x_{n+1-j}^{(n)}=0\quad\text{for}\ j=1,\ldots ,n.
\end{equation}
It is manifest from the explicit differential equations that the gradient system in Eqs. \eqref{g-flow}, \eqref{g-flow-chahn} preserves this parity symmetry. More specifically, if the initial condition satisfies
\begin{subequations}
\begin{equation}
x_j(0)+x_{n+1-j}(0)=0\quad\text{for}\ j=1,\ldots ,n,
\end{equation}
then so does the corresponding gradient flow $\bigl(x_1(t),\ldots ,x_n(t)\bigr)$, $t\geq 0$:
\begin{equation}
x_j(t)+x_{n+1-j}(t)=0\quad\text{for}\ j=1,\ldots ,n.
\end{equation}
\end{subequations}
If we now perform the reduction of our gradient system to the pertinent $m:=\lfloor \frac{n}{2}\rfloor$-dimensional invariant manifold through the substitution
\begin{equation}
(x_1,\ldots ,x_n)=\begin{cases}
(-y_m,\ldots ,-y_1,y_1,\ldots ,y_m)&\text{if}\ n=2m,\\
(-y_m,\ldots ,-y_1,0,y_1,\ldots ,y_m)&\text{if}\ n=2m+1,
\end{cases}
\end{equation}
then we arrive at the differential equations
\begin{subequations}
\begin{align}\label{gflow-symmetric-even}
{ \frac{\text{d} y_j}{\text{d} t}-\pi ( j -{\textstyle \frac{1}{2}})} &+
{ \arctan\left(\frac{y_j}{a}\right)  +\arctan\left(\frac{y_j}{b}\right) +\arctan\left( 2y_j \right) } \\
&+
 \sum_{\substack{1\leq k\leq m \\ k \neq j}} \bigl( \arctan (y_j+y_k)+\arctan (y_j-y_k) \bigr) =0 
 \nonumber
\end{align}
($j=1,\ldots ,m$) when $n=2m$, and
\begin{align}\label{gflow-symmetric-odd}
{ \frac{\text{d} y_j}{\text{d} t}-\pi j } &+
{ \arctan\left(\frac{y_j}{a}\right)  +\arctan\left(\frac{y_j}{b}\right) +\arctan\left( 2y_j \right)+\arctan\left( y_j\right) } \\
&+
 \sum_{\substack{1\leq k\leq m \\ k \neq j}} \bigl( \arctan (y_j+y_k)+\arctan (y_j-y_k) \bigr) =0 
 \nonumber
\end{align}
\end{subequations}
($j=1,\ldots ,m$) when $n=2m+1$, respectively. By comparing the latter differential equations with the gradient system for the Wilson polynomials in Eq. \eqref{g-flow-wilson}, it is seen that
for $n=2m$ we recover the case $(a,b,c,d)\to (a,b,\frac{1}{2},0)$ and for $n=2m+1$ we recover the case $(a,b,c,d)\to (a,b,\frac{1}{2},1)$. Indeed, the symmetric continuous Hahn polynomials and the Wilson polynomials are known to be related by the following quadratic relations (cf. e.g. \cite[Section 2.4]{koo:quadratic}):
\begin{equation}
\text{p}_{n}(x ;a,b)= \begin{cases} \text{p}_m(x^2 ;a,b,\frac{1}{2},0)&\text{if}\ n=2m,\\
x \text{p}_m(x^2 ;a,b,\frac{1}{2},1)&\text{if}\ n=2m+1.
\end{cases}
\end{equation}

The upshot is  that for initial conditions respecting the parity invariance,  the estimate for the rate of the exponential convergence in the second part of Theorem \ref{cHzeros:thm} can be improved as follows.

\begin{theorem}\label{cHzeros-symmetric:thm}
\begin{subequations}
Let $\text{Re}(a),\text{Re}(b)>0$ with possible non-real parameters  $a$ and $b$ arising as a complex conjugate pair, and
 let
$\bigl(x_1(t),\ldots ,x_n(t)\bigr)$, $t\geq 0$ denote the unique solution of the gradient system \eqref{g-flow}, \eqref{g-flow-chahn} determined by a choice for
the initial condition $(x_1(0),\ldots ,x_n(0))\in\mathbb{R}^n$ such that
\begin{equation*}
x_j(0)+x_{n+1-j}(0)=0\quad\text{for}\ j=1,\ldots ,n.
\end{equation*}
Then for any $\kappa$ in the interval
\begin{align}\label{cH-rate-symmetric}
{\textstyle 0< \kappa <   \frac{2 \lfloor \frac{n}{2}  \rfloor }{1+4R_n^2}+  \frac{\text{Re}(a)}{\text{Re}^2(a)+(R_n+|\text{Im}(a)|)^2}+ \frac{\text{Re}(b)}{\text{Re}^2(b)+(R_n+|\text{Im}(b)|)^2} } 
{\textstyle + \frac{1-(-1)^{n}}{2+2R_n^2}} ,
\end{align}
where $R_n:=x_n^{(n)}$, there exists a constant  $K>0$ (depending on the parameter values, on the initial condition, and on $\kappa$) such that 
\begin{equation}
|x_j(t)-x_j^{(n)}|\leq  K e^{-\kappa t} \qquad (j=1,\ldots,n)
\end{equation}
whenever $t\geq 0$.
\end{subequations}
\end{theorem}

For $n=2m+1$, Theorem \ref{cHzeros-symmetric:thm} is  immediate from the reduced differential equation \eqref{gflow-symmetric-odd} and
 the second part of Theorem \ref{Wzeros:thm}, upon performing the parameter specialization $(a,b,c,d)\to (a,b,\frac{1}{2},1)$. The case $n=2m$ would  follow similarly via
 the reduced differential equation \eqref{gflow-symmetric-even}  
 upon performing the parameter specialization $(a,b,c,d)\to (a,b,\frac{1}{2},0)$ in Theorem  \ref{Wzeros:thm}. However, since the parameter specialization $d\to 0$ takes us outside the Wilson parameter domain
 considered here,  for the argument to stick rigorously one formally has to repeat the proof  of Theorem  \ref{Wzeros:thm} for the case of the gradient flow in Eq. \eqref{gflow-symmetric-even}. 
 
\begin{figure}[h]
 \centering
 \caption{Continuous Hahn trajectories $x^{(30)}_1(t),\ldots, x_{30}^{(30)}(t)$  corresponding to the initial condition $x_j^{(30)}(0)=3$ ($j=1,\ldots ,30$) and the parameter values $a=10$ and $b=\frac{3}{10}$.}
\includegraphics[scale=0.7]{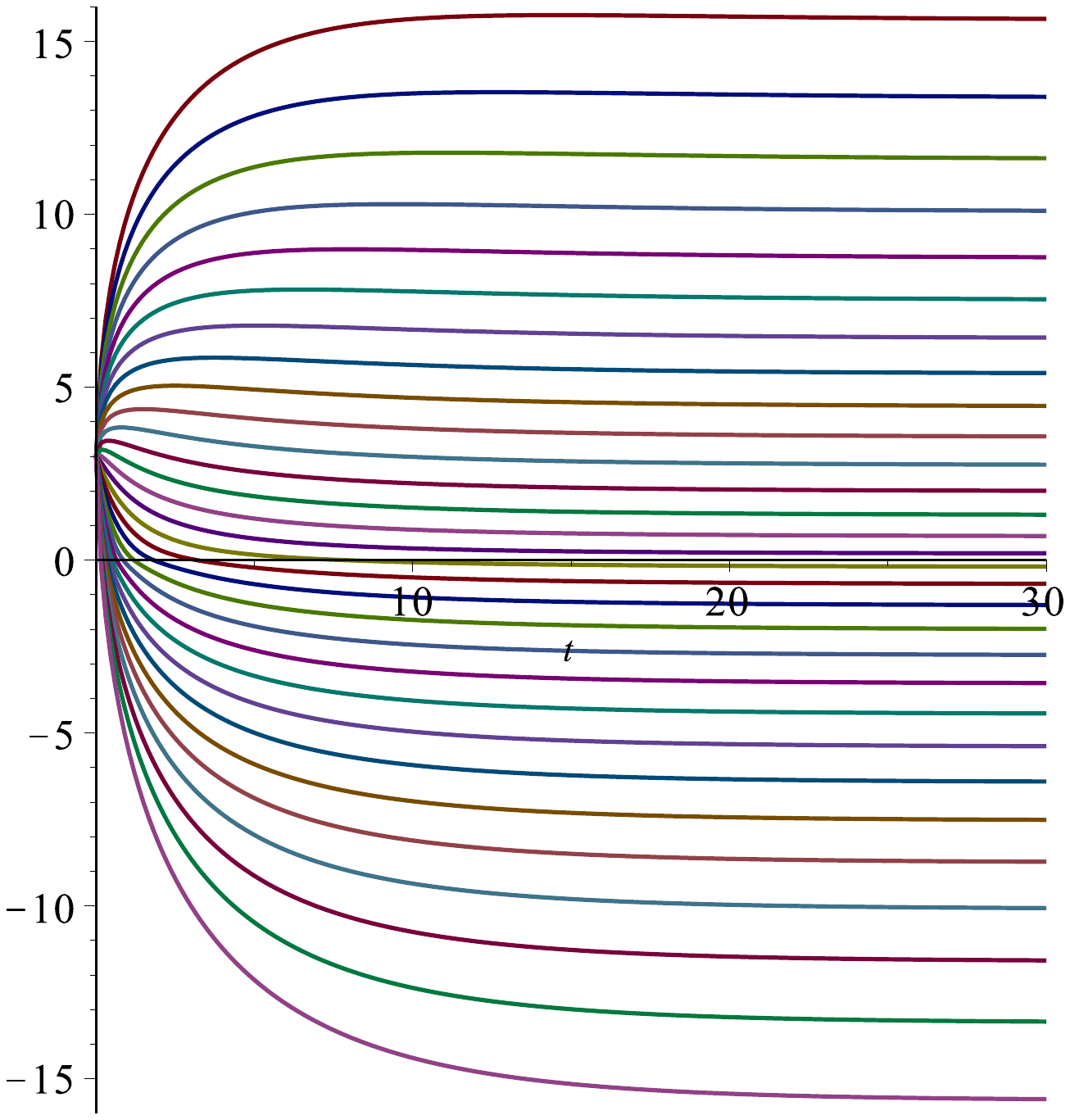} \label{cH-symmetric-fig1}
\end{figure}

\begin{figure}[h]
  \centering
  \caption{Evolution of the continuous Hahn logarithmic error  $\log\bigl |x_j(t)-x_j^{(30)}\bigr| $ for $j=1$ (top), $j=8$ (second from below), $j=23$ (bottom), and $j=30$ (second from above), where the ordering refers to that of the asymptotic tails, with the initial condition $x_j^{(30)}(0)=3$ ($j=1,\ldots ,30$) and the parameter values $a=10$ and $b=\frac{3}{10}$.}
\includegraphics[scale=0.7]{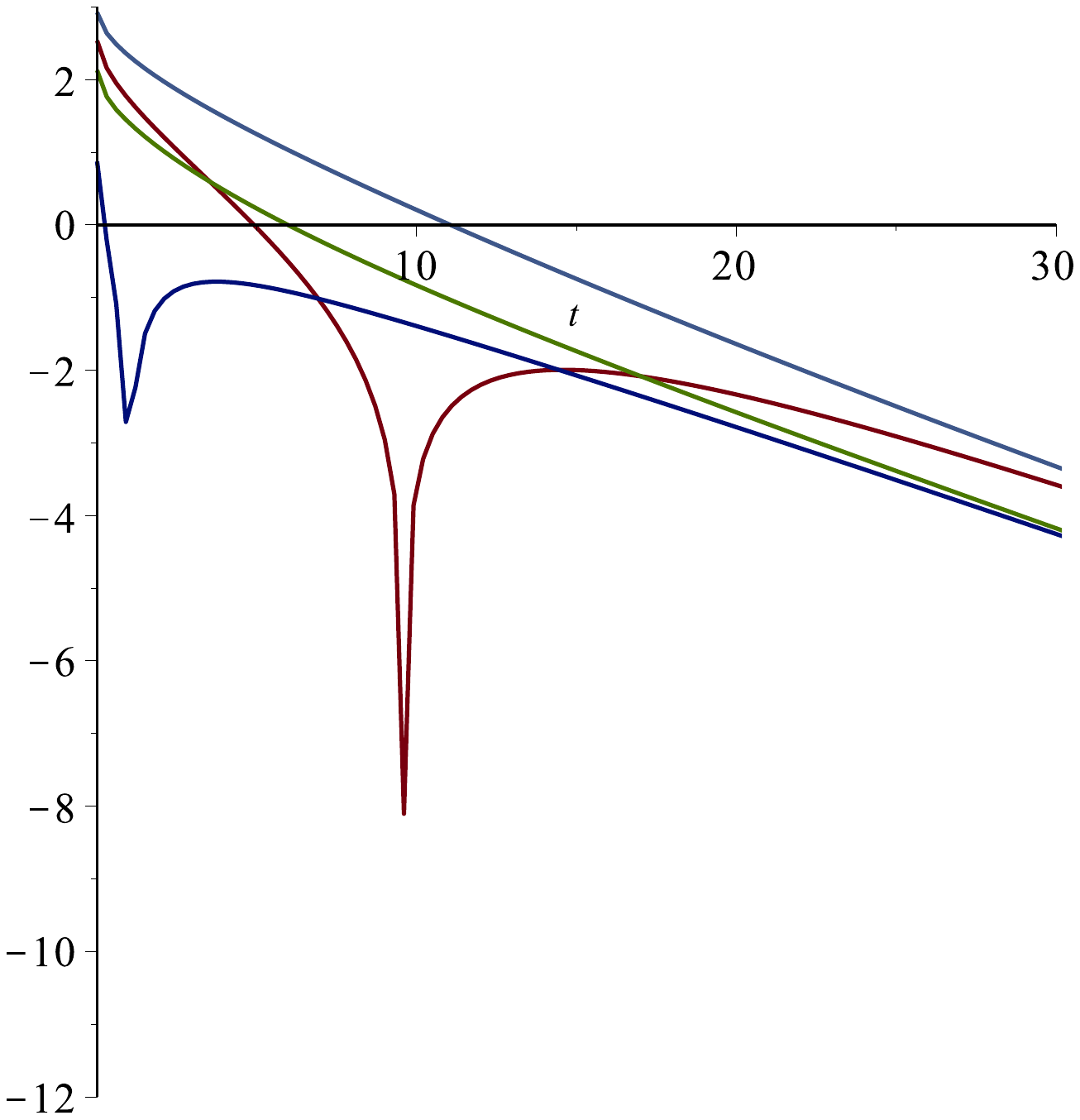}  \label{cH-symmetric-fig2}
\end{figure}

Since our numerical example in Section \ref{cH-numerics} corresponded to an initial condition of the form $x_j(0)=0$ ($j=1,\ldots ,n$), the improved  exponential convergence of Theorem \ref{cHzeros-symmetric:thm}  actually applies in this situation.  However, if we break the parity symmetry by moving the initial condition up to $x_j(0)=3$ ($j=1,\ldots ,n$)---while maintaining the parameter values $a=10$ and $b=\frac{3}{10}$
(cf. Figure \ref{cH-symmetric-fig1})---then we see from the slopes of the logarithmic error (cf. Figure \ref{cH-symmetric-fig2}) that
the rate of the exponential convergence indeed slows down considerably. Notice at this point that the downward peaks in the evolution of the logarithmic error reveal that the corresponding limiting values are no longer approached monotonously: the peak detects when the trajectory of $x_j(t)$ overshoots its limiting value $x_j^{(n)}$.

\section*{Acknowledgments}
It is a pleasure to thank Alexei Zhedanov for emphasizing that the Morse functions from  \cite{die-ems:solutions}, which minimize at the roots of the continuous Hahn, Wilson and Askey-Wilson polynomials, 
should be viewed as natural analogs of Stieltjes' electrostatic potentials for the roots of the classical orthogonal polynomials. 
Thanks are also due to an anonymous referee for suggesting some important improvements in the presentation.

\end{document}